\documentstyle[12pt,amscd]{amsart}
\newtheorem{thm}[equation]{Theorem}
\numberwithin{equation}{section}
\newtheorem{cor}[equation]{Corollary}

\newtheorem{lem}[equation]{Lemma}

\newtheorem{prop}[equation]{Proposition}

\begin{document}
\raggedbottom
\voffset=-.7truein
\hoffset=0truein
\vsize=8truein
\hsize=6truein
\textheight=8truein
\textwidth=6truein
\baselineskip=18truept
\def\mapright#1{\smash{\mathop{\longrightarrow}\limits^{#1}}}
\def\mapleft#1{\smash{\mathop{\longleftarrow}\limits^{#1}}}
\def\mapup#1{\Big\uparrow\rlap{$\vcenter {\hbox {$#1$}}$}}
\def\mapdown#1{\Big\downarrow\rlap{$\vcenter {\hbox {$\ssize{#1}$}}$}}
\def\mapne#1{\nearrow\rlap{$\vcenter {\hbox {$#1$}}$}}
\def\mapse#1{\searrow\rlap{$\vcenter {\hbox {$\ssize{#1}$}}$}}
\def\mapr#1{\smash{\mathop{\rightarrow}\limits^{#1}}}
\def\ss{\smallskip}
\def\sm{\wedge}
\def\la{\langle}
\def\ra{\rangle}
\def\on{\operatorname}
\def\kbar{{\overline k}}
\def\qed{\quad\rule{8pt}{8pt}\bigskip}
\def\ssize{\scriptstyle}
\def\a{\alpha}
\def\bz{{\bold Z}}
\def\vp{v_1^{-1}\pi}
\def\im{\on{im}}
\def\coef{\on{coef}}
\def\ext{\on{Ext}}
\def\bspin{\on{BSpin}}
\def\sq{\on{Sq}}
\def\eps{\epsilon}
\def\ar#1{\stackrel {#1}{\rightarrow}}
\def\bq{{\bold Q}}
\def\bk{{\bold K}}
\def\br{{\bold R}}
\def\bc{{\bold C}}
\def\si{\sigma}
\def\Ebar{{\overline E}}
\def\Sum{\sum}
\def\tfrac{\textstyle\frac}
\def\tb{\textstyle\binom}
\def\Si{\Sigma}
\def\w{\wedge}
\def\equ{\begin{equation}}
\def\b{\beta}
\def\G{\Gamma}
\def\g{\gamma}
\def\endeq{\end{equation}}
\def\sn{S^{2n+1}}
\def\zp{\bold Z_p}
\def\P{{\cal P}}
\def\zt{\bold Z_2}
\def\Hom{\on{Hom}}
\def\ker{\on{ker}}
\def\coker{\on{coker}}
\def\da{\downarrow}
\def\io{\iota}
\def\Om{\Omega}
\def\car{{\cal R}}
\def\u{{\cal U}}
\def\e{{\cal E}}
\def\exp{\on{exp}}
\def\xbar{{\overline x}}
\def\ebar{{\overline e}}
\def\et{{\widetilde E}}
\def\ni{\noindent}
\def\coef{\on{coef}}
\def\den{\on{den}}
\def\lcm{\on{l.c.m.}}
\def\vi{v_1^{-1}}
\def\ot{\otimes}
\def\psibar{{\overline\psi}}
\def\mhat{{\hat m}}
\def\exc{\on{exc}}
\def\ms{\medskip}
\def\ehat{{\hat e}}
\def\BSpt{\widetilde{BSp}}
\def\vt{{\widetilde v}}
\def\tbinom#1#2{{\textstyle\binom {#1}{#2}}}
\def\dirlim{\on{dirlim}}
\def\x^{X_K^{\wedge}}
\title[1-line for Spin(2n+1)]
{The 1-line of the $K$-theory Bousfield-Kan spectral sequence for $Spin(2n+1)$}
\author[Bendersky]{Martin Bendersky}
\address{Hunter College,
CUNY, NY, NY 10021}
\email{mbenders@@shiva.hunter.cuny.edu}
\author[Davis]{Donald M. Davis}
\address{Lehigh University\\Bethlehem, PA 18015}
\email{dmd1@@lehigh.edu}

\date{September 20, 1999}

\maketitle

\section{Statement of results} \label{intro}
The $p$-primary $v_1$-periodic homotopy groups, $\vp_*(X;p)$, of a topological
space $X$, as defined in \cite{DM}, are a localization of the portion of the
actual homotopy groups of $X$ detected by $K$-theory. In \cite{D,BD2,BDM},
$\vp_*(X;p)$ was calculated for classical groups $X$ and primes $p$ in all
cases except $(SO(n),2)$.  In this paper we make a first step toward the
calculation of $\vp_*(SO(n);2)$, which is of course isomorphic to
$\vp_*(Spin(n);2)$.

In \cite{BT}, a Bousfield-Kan-type spectral sequence $E_r^{s,t}(X)$ based on
periodic $K$-theory was introduced. It converges to
the homotopy groups of the $K$-completion $\x^$ of the space $X$.
In this paper, we will deal exclusively with the localization of this spectral
sequence at the prime 2. We expect to prove in a subsequent paper that if
$X=Spin(n)$ then
$\pi_i(\x^)\approx\vp_i(X;2)$ for sufficiently large values of $i$.
Thus, since $v_1$-periodic homotopy groups are periodic, a computation
of this spectral sequence for $X=Spin(n)$ and $t-s$ large
would yield a complete computation of $\vp_*(Spin(n))$.

  If $X$ is an $H$-space such that $K_*(X)$ is a
free commutative algebra on odd-dimensional classes (e.g. $X=Spin(n)$), then
(\cite[4.9]{BT})
 $E_2(X)$ is the homology of an unstable cobar complex
determined by these classes and their $K_*K$-coaction. The main result
of this paper is an explicit computation of the 1-line $E_2^{1,t}(Spin(2n+1))$
of this spectral
sequence when $p=2$ and $X=Spin(2n+1)$, .

The first step toward this computation is the following general result,
which was inspired by the results of Bousfield (\cite{Bo}) at the odd primes.
This result will be proved in Section \ref{Thm1}. In it, $(-)^\#$ denotes the
Pontryagin dual, and $PK^1(X)$ the primitives in $K^1(X)$.
\begin{thm} \label{E2iso}If $X$ is a simply-connected finite $H$-space with
$K_*X$ an exterior algebra on odd-dimensional classes of degree $\le2M+1$,
then
$E_2^{1,2m}(X)=0$ for all $m$, and, if $m>M$, then
$$E_2^{1,2m+1}(X)\approx (PK^1(X)/\im(\psi^r-r^m:r=2,3,-1))^\#.$$
\end{thm}

For simplicity of notation, we define
$$\vt^{2m}(X):=PK^1(X)/\im(\psi^r-r^m:r=2,3,-1)$$
and the closely-related functor
$$v^{2m}(X):=PK^1(X)/\im(\psi^2,\psi^3-3^m,\psi^{-1}-(-1)^m).$$
The advantage of $\vt^*(-)$ is that it is more closely related to
$E_2^{1,*+1}(-)$, while the advantage of $v^*(-)$ is that it is periodic
and has simpler formulas. They are related by the following result, which will
be proved at the end of this section.
\begin{prop}\label{vv} If $X$ is a simply-connected finite $H$-space with
$K_*X$ an exterior algebra on classes of degrees $2d_1+1,\ldots,2d_r+1$,
then $v^{2m}(X)\approx \vt^{2m}(X)$ if $m>\sum d_i$.
\end{prop}

 The second purpose of this
note is to compute $v^{2m}(Spin(2n+1))$, which, by Theorem \ref{E2iso} and
Proposition \ref{vv},
is isomorphic to $E_2^{1,2m+1}(Spin(2n+1))$ if $m>n^2$. Here we use the
rational equivalence of $Spin(2n+1)$ with $\prod_{i=1}^n S^{4i-1}$, and that
the sum of the odd integers up to $2n-1$ is $n^2$.

We will show in Proposition \ref{Bmor}
that there is an injection of Adams-modules
$$PK^1(Sp(n))\to PK^1(Spin(2n+1)),$$
and hence a morphism
\begin{equation}\label{vmor}v^{2m}
(Sp(n))\to v^{2m}(Spin(2n+1)),\end{equation}
which is dual to a morphism of $E_2^{1,2m+1}(-)$. The following result was
proved in \cite{BDspin}. Here and throughout, $\nu(-)$ denotes the exponent of
2 in a number.
\begin{thm}$($\cite{BDspin}$)$ If $m$ is odd and $m\ge 2n$, then
$E_2^{1,2m+1}(Sp(n))\approx\bz/2^{eSp(m,n)}$, where
$$eSp(m,n):=\min(\nu(S_{m,j}):j>2n),$$
with $S_{m,j}=\sum(-1)^k\tbinom jkk^m$.
If $m$ is even, then $E_2^{1,2m+1}(Sp(n))=\bz/2$.\label{Spthm}
\end{thm}
\noindent Actually, the result in \cite{BDspin} was proved for the $BP$-based
unstable Novikov spectral sequence,
but one might use the change of rings theorem
in \cite{BT}, or  mimic the calculation in \cite{BDspin} to
yield the result for the $K$-based spectral sequence.

Our second main result is as follows. In light of \ref{E2iso}, it gives
$E_2^{1,2m+1}(Spin(2n+1))$ when $m>n^2$.
\begin{thm}\label{main} If $n=3$ or $n\ge5$, and $m$ is odd, then
$$v^{2m}(Spin(2n+1))\approx\bz/2^{\min(eSp(m,n),\nu(R_1(m,n)),\nu(R_2(m,n)))}
\oplus\bz/2^{\min(2+\nu(m+1),n)},$$
where
$$R_1(m,n)=\sum_{\text{odd }k\ge1}k^m\biggl(\sum_{i=0}^{n-k}\tbinom{2n+1}i-
2\sum_{t\ge0}\tbinom{2n+2}{n-1-k-4t}\biggr)$$
and $R_2(m,n)=$
$$\bigg((2^{2n+1}-3^{m+1}+1)\sum_{\text{odd }k\ge1}k^m\sum_{t\ge0}
\tbinom{2n+2}{n-1-k-4t}-3\cdot2^{2n}\sum_{\text{odd }k\ge1}k^m\sum_{t\ge0}
\tbinom{2n+1}{n-1-k-3t}\biggr)\biggm/2^{\min(2+\nu(m+1),n)}.$$
If $m$ is even, then $v^{2m}(Spin(2n+1))\approx\bz/2\oplus\bz/2$.
\end{thm}
The sums over odd $k$
could include even values as well, as long as $m>n^2$, since this is
greater than the largest possible exponent of $v^*(Spin(2n+1))$.
We will discuss the slight anomaly when $n=4$ in
Proposition \ref{n=4}. Since $Spin(5)\approx Sp(2)$, the case $n=2$ has been
covered in \cite{BDM} (and is quite different).

The morphism (\ref{vmor}) sends a generator to $g_1+ag_2$, for some integer
$a$. We could certainly determine $a$, but shall not bother to do so here.

We illustrate with concrete calculations when $n=5$ and $6$. These results are
obtained by computer calculation, although we shall show after (\ref{bigsum})
how to obtain
some of them by hand.
\begin{prop} For $n=5$ and $6$, $v^{2m}(Spin(2n+1))\approx\bz^{e_1(2n+1)}\oplus
\bz^{e_2(2n+1)}$, where $e_1(2n+1)$ and $e_2(2n+1)$ are as in the following
table, which also presents $eSp(m,n)$ for comparison.\label{5and6}\end{prop}

\begin{tabular}
{|cc|ccc|}
\hline
$n$&$m$&$eSp(m,n)$&$e_1(2n+1)$&$e_2(2n+1)$\\ \hline
$5$&$3\mod 8$&$8$&$5$&4\\
5&$7 \mod 8$&$\min(11,\nu(m-7)+6)$&$\min(11,\nu(m-103)+2)$&5\\
5&$1 \mod 4$&$\min(14,\nu(m-73)+6)$&$\min(12,\nu(m-73)+4)$&3\\ \hline
6&$15 \mod 16$&11&6&6\\
6&$7 \mod 16$&11&$\min(9,\nu(m-23)+4)$&5\\
6&$3 \mod 8$&$\min(15,\nu(m-75)+9)$&$\min(15,\nu(m-523)+5)$&4\\
6&$1 \mod 4$&$\min(14,\nu(m-9)+8)$&$\min(13,\nu(m-201)+5)$&3\\
\hline
\end{tabular}

\begin{pf*}{Proof of Proposition \ref{vv}} There is a rational equivalence
$X\simeq\prod S^{2d_i+1}$. Thus there are elements $x_i\in PK^1(X)$ for which
$\psi^k(x_i)=k^{d_i}x_i$ for all $k$. Similarly to the procedure in
\cite{Drep}, we can find a basis $\{w_1,\ldots,w_r\}$ of $PK^1(X)_{(2)}$
such that each $w_i$ is a rational linear combination of $\{x_j:j\ge i\}$, and
hence the matrix of each $\psi^k$ is triangular with
$\{k^{d_1},\ldots,k^{d_r}\}$ along the diagonal. Thus the Adams module
$PK^1(X)$ can be built from $PK^1(S^{2d_1+1}),\ldots,PK^1(S^{2d_r+1})$
by short exact sequences, and by the Snake Lemma this implies the same of
$v^{2m}(X)$ (built from $v^{2m}(S^{2d_i+1})$). Since
$2^{d_i}v^{2m}(S^{2d_i+1})=0$ by \cite{BT}, we deduce $2^{\sum d_i}v^{2m}(X)=0$
(a very conservative estimate).
If $m\ge\sum d_i$, then $\cdot2^m$ is 0 in $v^{2m}(X)$. Since the only
difference between $v^{2m}(X)$ and $\vt^{2m}(X)$ is multiples of $2^m$, we
conclude that they are equal.\end{pf*}

\section{Proof of Theorem \ref{E2iso}}\label{Thm1}
In this section, we prove Theorem \ref{E2iso}.
We begin by observing that the exact sequence in $E_2$ induced by the
coefficient sequence $$0\to \bz\to \bq\to \bq/\bz\to 0$$
induces isomorphisms $E_2^{s,t}(X;\bq/\bz)\to E_2^{s+1,t}(X)$
except when $s=0$, $t=2d_i+1$. Hence there is an isomorphism when
 $s=0$, $t>2d_r+1$. Everything is localized at 2. We have
\begin{equation}\label{K}E_2^{0,2m+1}(X;\bq/\bz)\approx
\ker\bigl(\bar d: K_{2m+1}(X)\ot\bq/\bz\to (K_*K\ot K_*X\ot\bq)/
 U(X)\bigr),\end{equation}
where $\bar d$ is induced from
$$d:(K_*(X)\ot\bq,K_*(X))\to (K_*K\ot K_*X\ot\bq,U(X)),$$  the
boundary
in the unstable cobar complex described in \cite[\S4]{BT}. The following
result, whose proof is deferred until the end of this section,
is crucial. Here $QK_*(-)$ denotes the indecomposables.
\begin{prop} Suppose $QK_*(X)$ has basis $\{y_1,\ldots,y_r\}$ as a
$K_*$-module with $y_j\in K_{2d_j+1}(X)$, and $x\in K_*(X)\ot\bq$ has
$d(x)=\sum s_j\ot y_j$, with
$s_j\in K_*K\ot\bq$.
Then $x\in \ker({\bar d})$ $($of $(\ref{K}))$ if and only if for all positive
integers $k$ and all subscripts $j$
\begin{equation}k^{d_j}\langle \psi^k,s_j\rangle \in {K_*}_{(2)}.
\label{unscond}\end{equation}
\label{unstable}\end{prop}
The condition (\ref{unscond}) says that each $k^{d_j}\langle \psi^k,s_j\rangle$
is a 2-local integer times the generator of $K_{|s_j|}$.

Now we prove Theorem \ref{E2iso}.
Suppose $\{y_1,\ldots,y_r\}$ is as in Proposition \ref{unstable} and
$$\psi:K_*X\to K_*K\ot K_*X$$ satisfies
\begin{equation}\label{psi}\psi(y_j)=\sum_i s_{i,j}\ot y_i.\end{equation}
Let $u\in K_2$ denote the Bott class, so that $K_*=\bz_{(2)}[u]$, and
let $v=\eta_R(u)$.
Then $$d(u^ey_j)=(v^e-u^e)\ot y_j-\sum_{i\ne j}u^es_{i,j}\ot y_i.$$
Since, by \cite[p.676]{Boa}, $\langle \psi^k,v\rangle=ku$ and $\langle
\psi^k,u\rangle=u$, the expression whose coefficients correspond to
$\langle \psi^k,s_j\rangle$ in (\ref{unscond}) equals
$$k^eu^ey_j-\sum_{\text{all }i}\langle\psi^k,s_{i,j}\rangle u^ey_i.$$

A basis for $QK_{2m+1}(X)$ is given by $\{u^{m-d_1}y_1,\ldots,u^{m-d_r}y_r\}$,
and, by Proposition \ref{unstable}, $E_2^{0,2m+1}(X;\bq/\bz)$ is the
intersection, over all integers $k>1$,
 of the kernel of the morphism
$$\delta_k:K_{2m+1}(X)\ot\bq/\bz\to K_{2m+1}(X)\ot\bq/\bz$$
defined by
$$\delta_k(u^{m-d_j}y_j)=-k^{m-d_j}u^{m-d_j}k^{d_j}y_j+u^{m-d_j}\sum_i\langle
\psi^k,s_{i,j}\rangle k^{d_i}y_i.$$
Here we have negated the differential $d$ for later convenience. Define $\a_{i,
j,k}\in\bz$
by $\langle\psi^k,s_{i,j}\rangle k^{d_i}=\a_{i,j,k}u^{d_j-d_i}$, and define an
$r\times r$-matrix $A_k=(\a_{i,j,k})$. Then $\delta_k$ is a linear
transformation of free $\bq/\bz$-modules with matrix $A_k-k^mI$.

Now let $w_i\in K^{2d_i+1}(X)$ be dual to $y_i$, and let $z_i=u^{d_i}w_i\in
K^1(X)$. By \cite[11.19]{Boa} applied to (\ref{psi}), we have
$$\psi^k(w_i)=\sum_j \langle\psi^k,s_{i,j}\rangle w_j.$$
Thus
$$\psi^k(z_i)=k^{d_i}u^{d_i}\sum_j\langle\psi^k,s_{i,j}\rangle u^{-d_j}z_j
=
\sum_j\a_{i,j,k}z_j.$$
Thus $A_k-k^mI$ is the matrix whose $i$th row gives $(\psi^k-k^m)(z_i)$
expressed in terms  of the basis $\{z_1,\ldots,z_r\}$.
 Then
$\vt^{2m}(X)$ is the abelian group presented by a matrix obtained
by stacking the matrices $A_k-k^mI$ for all $k$.
Thus we see that the same matrix $A$, defined to be all matrices
$A_k-k^mI$ stacked, yields
$E_2^{0,2m+1}(X;\bq/\bz)$ if we take the kernel of the linear transformation of
free $\bq/\bz$-modules whose matrix is $A$, and it yields $\vt^{2m}(X)$ if we
think of it as the presentation of an abelian group. Theorem \ref{E2iso}
follows now from the following proposition.
\begin{prop} If $A$ is an integer matrix, let $K(A)$ denote the kernel of the
linear transformation of free $\bq/\bz$-modules whose matrix is $A$, and let
$G(A)$ denote the abelian group presented by $A$. Then there is an isomorphism
of abelian groups
$K(A)\approx (G(A))^{\#}$. The same is true if we are localized
at a prime $p$.\end{prop}
\begin{pf} The matrix $A$ can be brought to diagonal form (with rows and
columns of 0's adjoined) by integer row and
column operations. These operations consist of interchanging, multiplying by $-
1$, and adding a multiple of one to another. For the first interpretation of
the matrix, they correspond to change-of-basis in the domain and range, and
hence do not affect the kernel (up to isomorphism),
while in the second interpretation they
correspond to invertible change of generating set and to simplification of
relations, and hence do not affect the isomorphism class of the group
presented. Let $D$ be a
diagonal matrix with entries $d_i>1$, $m$ 1's, and $n$ 0's on the diagonal, and
possibly many additional rows of 0's. Then
$K(D)\approx
\bigoplus \bz/d_i\oplus (\bq/\bz)^n$, while $G(D)\approx\bigoplus\bz/d_i
\oplus (\bz)^n$. Since $(\bz)^{\#}\approx\bq/\bz$ and $(\bz/d)^{\#}\approx
\bz/d$, the proposition follows.\end{pf}

\begin{pf*}{Proof of Proposition \ref{unstable}} We must show that if
$d(x)=\sum s_j\ot y_j$, then ${\bar d}(x)\in U(X)$ if and only if
$k^{d_j}\langle \psi^k,s_j\rangle$ is 2-integral for all $j$.

The following analysis of the unstable condition (i.e. the condition for being
in $U(X)$) is based on \cite[\S7]{BT}, which derived from \cite{BCRS}.
Let $\bk_{2n}$ denote the $2n$th space in the $\Omega$-spectrum for $K$.
There is a homotopy equivalence $\bk_{2n}\to BU$, and the structure maps
$\Sigma^2\bk_{2n}\to\bk_{2n+2}$
are the Bott maps $B$. These combine to a tower
$$K_{2e}(BU)\mapright{B_*}
K_{2e+2}(BU)\mapright{B_*}
\cdots\mapright{B_*}K_{2e+2d_j}(BU)\mapright{B_*}\cdots\to
K_{2e}(K),$$
with $K_{2e}(K)$ the direct limit of the tower.
There is a similar tower after tensoring with $\bq$.
Suppose $s_j\in K_{2e}(K)\ot\bq$.
Then $s_j\ot y_j\in U(X)$ if and only
if $s_j$ pulls back to an element of $K_{2e+2d_j}(BU)$. This follows from the
description of unstable comodules in the fourth paragraph after \cite[4.3]{BT}.

For purposes of calculation, it is convenient to write the above tower in the
equivalent form
$$K_0(BU)\mapright{B_*}K_0(BU)\mapright{B_*}\cdots\mapright{B_*}K_0(BU)
\mapright{B_*}\cdots\to K_0(K).$$
In \cite[1.2-1.4]{BCRS}, it is
observed that $K_0(BU) \subset \bq[w]$ is a polynomial algebra
generated by  $\binom wn =w(w-1)\cdots (w-n+1)/n!\in \bq[w]$ for $n\ge1$. Let $A
$ denote the free abelian group on
the rational polynomials  $\binom wn$. Then
$A$ consists of all rational polynomials $f(w)$
such that $f(n)$ is an integer for every integer $n$, so-called {\em numerical
polynomials}.
Then, still following \cite{BCRS},
$B_*(f)-wf$ is decomposable for all $f\in A$,
and $B_*$ annihilates decomposables, from which it is  deduced that
$K_0(K)\approx w^{-1}A$, with an element $f$ in the $d$th factor corresponding
to $w^{-d}f$ in the limit.

It is also noted in \cite[p.390]{BCRS} that $w=u^{-1}v\in K_0(K)$. This implies
that $\langle\psi^k,w^n\rangle=k^n$. Thus $\langle\psi^k,f(w)\rangle=f(k)$,
and so $A$ consists of those polynomials $f(w)$ for which $\langle\psi^k,f(w)
\rangle$ is integral for all $k$.

Now $s_j\in K_0(K)\ot\bq$ is represented by a finite Laurent series $s_j(w)$ and
pulls back to an element in the $d_j$th copy of $K_0(BU)$ if and only if
$w^{d_j}s_j(w)$ is a numerical polynomial if and only if
$\langle\psi^k,w^{d_j}s_j(w)\rangle$ is integral for all $k$. Since
$\langle\psi^k,w^{d_j}s_j\rangle=k^{d_j}\langle\psi^k,s_j\rangle$, our desired
conclusion follows.
\end{pf*}

\section{Proof of Theorem \ref{main}}\label{proof}
We begin with a result culled from \cite{Nay}.
\begin{prop}\label{Bmor} There are morphisms of Adams modules
\begin{eqnarray*} PK^1(SU(2n+1))\mapright{j^*}&PK^1(Sp(n))&\mapright{\theta}
PK^1(Spin(2n+1))\\
\langle B_1,\ldots,B_{2n}\rangle\longrightarrow&\langle B_1',\ldots,B_n'\rangle
&\longrightarrow\langle B_1'',\ldots,B_{n-1}'',D\rangle,
\end{eqnarray*}
where $\langle-,\ldots,-\rangle$ denotes the free abelian group on the
indicated classes,  $\theta(B_i')=B_i''$ if $i<n$,
$\theta(B_1'+\cdots+B_n')=2^{n+1}D$, and
$$j^*(B_i)=\begin{cases}B_i'&i\le n\\B'_{2n+1-
i}&i>n.\end{cases}$$
\end{prop}
\begin{pf} We use Hodgkin's result (\cite{Hod}) that $PK^1(G)\approx
\langle \b(\rho_1),\ldots,\b(\rho_l)\rangle$, where $\rho_1,\ldots,\rho_l$
are the fundamental representations of $G$, and $\b(\rho)$ is the
virtual vector bundle over $\Sigma G$ associated to $\rho-\dim(\rho)$.
Naylor (\cite[p.151]{Nay}) notes that the maps
$$Sp(n)\mapright{k}SU(2n)\mapright{i}SU(2n+1)$$
induce morphisms of representation rings
$$R(SU(2n+1))\mapright{i^*}R(SU(2n))\mapright{k^*}R(Sp(n))$$
satisfying $$k^*(\mu_i')=\begin{cases}\eta_i&i\le n\\
\eta_{2n-i}&i>n\end{cases}$$ and $i^*(\mu_i)=\mu_i'+\mu_{i-1}'$, where
$\eta_i$, $1\le i\le n$, $\mu_i'$, $1\le i\le 2n-1$,
and $\mu_i$, $1\le i\le
2n$, are the representations given by exterior power operations on
a canonical representation. We have $\mu'_0=\mu'_{2n}=1$.
Letting $B_i'=\b(\eta_i)+\b(\eta_{i-1})$
and $B_i=\b(\mu_i)$
yields the desired result about $j^*$, which is the composite $k^*\circ i^*$.

Naylor also considers the composite
$$j':Spin(2n+1)\to SO(2n+1)\to SU(2n+1).$$
We let $\lambda_i$ denote the exterior powers of the canonical unitary
representation of $Spin(2n+1)$, $\Delta$ the Spin representation,
$B_i''=\b(\lambda_i)$, and $D=\b(\Delta)$. Naylor shows that
$${j'}^*(B_i)=\begin{cases}B_i''&i\le n\\
B_{2n+1-i}''&i>n\end{cases}$$
and $2^{n+1}D=B''_1+\cdots+B''_n$.
Since $\ker(j^*)=\ker({j'}^*)$,
there is a homomorphism of abelian groups
$$PK^1(Sp(n))\mapright{\theta}PK^1(Spin(2n+1))$$
defined by $\theta(B_i')=B_i''$ for $i\le n$,
 and it will be an Adams-module homomorphism
since $j^*$ and ${j'}^*$ are Adams-module homomorphisms with $j^*$
surjective.
\end{pf}

Next we relate the above basis of $PK^1(SU(n))$ with one used by
Bousfield (\cite{Bo}).
\begin{prop}\label{bases} There is an isomorphism of Adams modules
\begin{equation}\label{CPiso}PK^1(SU(n))\approx \widetilde K^0(CP^{n-1}).
\end{equation}
Let $\xi_k$ denote the element of $PK^1(SU(n))$ which corresponds to $\xi^k-1$
under this isomorphism. Here $\xi$ denotes the Hopf bundle over $CP^{n-1}$.
Then
$$B_j=\sum(-1)^{k+1}\tbinom n{j-k}\xi_k.$$
\end{prop}
\begin{pf} If $I$ denotes the augmentation ideal of $R(G)$, then there
is an isomorphism $I/I^2\approx PK^1(G)$ under which $(-1)^{k+1}\lambda^k$
corresponds to $\psi^k$. (See, e.g., \cite[2.1ff]{Drep}.) We first claim that
under this isomorphism and (\ref{CPiso}) (which is well-known)
$\{\theta_n-n\}\in I/I^2$ corresponds to $\{\xi-1\}\in\widetilde K^0(CP^{n-
1})$. To see this, we use a result in \cite[p.206]{MT} which states that there
is a commutative diagram
$$\begin{CD} CP^{n-1}@>f>> BU\\
@VjVV @VVBV\\
\Om SU(n)@>\Om i>>\Om SU
\end{CD}$$
in which $j$ is a canonical map, $i$ is the inclusion, $B$ the Bott map, and
$f$ satisfies $f^*(c_i)=0$ unless $i=1$, and $f^*(c_1)$ is the usual generator
of $H^2(CP^{n-1})$.
This implies that $f$ classifies $\xi$. The isomorphism (\ref{CPiso})
is defined using the maps $j$ and $B$. Hence $\xi$ corresponds under this
isomorphism to the map
$i$, which is the canonical representation $\theta_n$.

Since $\psi^k(\xi)=\xi^k$, we obtain that
$$\xi_k=(-1)^{k+1}\lambda^k(\theta_n-n)=(-1)^{k+1}\sum_{j=0}^k\tbinom{-n}{k-j}
\lambda^j(\theta_n).$$
But $\lambda^j(\theta_n)=B_j$. To invert this equation, we take a formal sum of
the equations for all values of $k$, using powers of
an indeterminate $x$. We have
$$\sum(-1)^{k+1}\xi_k x^k=\sum x^k\sum_{j=0}^k\tbinom{-n}{k-j}B_j
=(1+x)^{-n}\sum x^jB_j$$
and hence
$$\sum x^jB_j=\sum_i\tbinom ni x^i\sum_k(-1)^{k+1}\xi_kx^k=\sum
x^j\sum_k\tbinom n{j-k}(-1)^{k+1}\xi_k,$$
as desired.
\end{pf}

The following corollary is immediate from the description of  $j^*$ in
Proposition \ref{Bmor} along with Proposition \ref{bases}.
\begin{cor}\label{Spn}
$PK^1(Sp(n))\approx\la\xi_1,\xi_2,
\ldots:R_{n+1},\ldots,R_{2n},S_j, j>2n\ra$, where
$$R_j=\sum (-1)^{k+1}\tbinom{2n+1}{j-k}\xi_k-\sum(-1)^{k+1}\tbinom{2n+1}{2n+1-
j-k}\xi_k,$$
$S_j=\sum(-1)^k\binom jk\xi_k$, and
$\psi^t\xi_k=\xi_{kt}$.
\end{cor}

Note that the relations allow one to express each $\xi_i$ with $i>n$ in terms 
of $\xi_1,\ldots,\xi_n$.
From Corollary \ref{Spn}, we easily deduce the following result.
\begin{cor}\label{vSpn} $v^{2m}(Sp(n))\approx\bz/2^{eSp'(m,n)}$, where
$$eSp'(m,n):=\min(
\nu(R_{m,n+1}),\ldots,\nu(R_{m,2n}),\nu(S'_{m,j}),j>2n),$$ with $R_{m,j}$ and
$S'_{m,j}$ defined by
$$R_{m,j}= \sum_{k \text{ odd}}\tbinom{2n+1}{j-k}k^m-\sum
_{k\text{ odd}}\tbinom{2n+1}{2n+1-j-
k}k^m,$$
$$S'_{m,j}= \sum_{k\text{ odd}}\tbinom jk k^m.$$
\end{cor}
\begin{pf} The sums are taken over all odd values of $k$ for which the binomial
coefficients are nonzero. The generator of the group is $\xi_1$, and the
relations are obtained from those in \ref{Spn}
by $\xi_k=\psi^2(\xi_{k/2})\equiv0$ if $k$ is even, while if $k$ is odd,
$\xi_k=\psi^k(\xi_1)\equiv k^m\xi_1$. Here we are using Corollary \ref{alls}
to imply that $\psi^k\equiv k^m$ in $v^{2m}(-)$ for all odd $k$.
\end{pf}

Because of Theorems \ref{E2iso} and \ref{Spthm} and Corollary \ref{vSpn}, 
we have
\begin{prop}\label{SUconj}
 $eSp'(m,n)=eSp(m,n)$; i.e.,
$$\min(
\nu(R_{m,n+1}),\ldots,\nu(R_{m,2n}),\nu(S_{m,j}),j>2n)=\min(
\nu(S_{m,j}),j>2n),$$ and so the relations $R_{m,j}$ are superfluous in $eSp'(m
,n)$
\end{prop}
We should remark on the minor difference between $S_{m,j}$ in \ref{Spthm} and
$S'_{m,j}$ in \ref{vSpn}, in that one involves a sum over all values of $k$,
while the other involves a sum only over odd values of $k$. As remarked in
Section \ref{intro}, we only consider values of $m$
large enough that summands which
are multiples of $2^m$ would not affect divisibility by 2 of the sum.
Thus $\nu(S_{m,j})=\nu(S'_{m,j})$ for such values of $m$.

Next we explain why the relations $\psi^3-3^m$ and $\psi^{-1}-(-1)^m$
imply relations
$\psi^k-k^m$ for all odd integers $k$. We contrast with the situation when
we are localized at an odd prime. The difference is due to the following
result, which is stated in \cite[2.9]{JX}.
\begin{prop}\label{units} If $p$ is an odd prime, then
the multiplicative group of units $(\bz/p^n)^\times$
is cyclic of order $(p-1)p^{n-1}$. If $r$ generates $(\bz/p^2)^\times$, then
any integer congruent to $r$ mod $p^2$ generates $(\bz/p^n)^\times$.
The group $(\bz/2^n)^\times$ is the direct product of the subgroup consisting of
$\pm1$ and the subgroup of classes congruent to 1 mod 4; the latter subgroup is
cyclic of order $2^{n-2}$ generated by any integer congruent to 5 mod 8.
\end{prop}

This implies that, for any $n$, any odd integer is congruent mod $p^n$ to
$\pm3^e$ for some $e$.
We also need the following result of Adams.
\begin{prop} $(\cite[5.1]{JXIII})$ If $X$ is a finite complex, there exists
an integer $e$ so that for all $k$ and $m$, $\psi^k\equiv \psi^{k+m^e}$ in
$K(X)/mK(X)$.\label{perio}
\end{prop}

We will use the following consequence of these results.
\begin{cor}\label{alls} Let $X$ be a finite complex, and if $p$ is odd let
$r$ generate $(\bz/p^2)^\times$. Let
$$Q_m=\begin{cases}PK^1(X;\bz_{(p)})/(\im(\psi^r-r^m))&\text{if $p$ is odd}\\
PK^1(X;\bz_{(2)})/(\im(\psi^3-3^m,\psi^{-1}-(-1)^m))&\text{if $p=2$.}\end{cases}
$$
Then $\psi^s=s^m$ in $Q_m$ for all $s\not\equiv 0$ mod $p$.
\end{cor}
\begin{pf} We consider just the case $p=2$. The case when $p$ is odd is similar
and slightly easier. We show $\psi^s\equiv s^m$ mod $2^n$ for all positive
integers $n$. This implies that they are equal in the finitely generated
$\bz_{(2)}$-module.

We first show that Proposition \ref{perio} is valid for $K^1(-)$ provided
$k$ and $m$ are relatively prime.
Adams' result dealt with $K^0(-)$. The validity in $K^{-1}(-)$ is immediate
since $K^{-1}(X)=K^0(\Sigma X)$. Since $\psi^k$ in $K^{-1}(-)$ corresponds to
$k\psi^k$ in $K^1(-)$, the periodicity in $K^{1}(-)$ follows, as long as
$k$ is a unit mod $m$.

Choose an integer $e$ that works for $X$ in Proposition \ref{perio}.
Use Proposition \ref{units} to find an integer $\ell$ and $\eps=0$ or 1
so that $s\equiv (-1)^\eps3^\ell$ mod $2^{ne}$. Then, for some $k$, we have,
in $Q_m$ mod $2^n$,
$$\psi^s=\psi^{(-1)^\eps3^\ell+k2^{ne}}\equiv \psi^{(-1)^\eps}\psi^{3^\ell}
\equiv ((-1)^\eps)^m(3^\ell)^m\equiv s^m.$$
\end{pf}

The following result is immediate from Propositions \ref{Bmor} and \ref{bases}.
\begin{prop}\label{PKSpin} The group
$PK^1(Spin(2n+1))$ has a subgroup isomorphic to
$PK^1(Sp(n))$, as described in Corollary \ref{Spn},
and an additional generator
$D$ satisfying
\begin{equation}\label{Ddef}
2^{n+1}D=\sum_{j=1}^n\sum_{k=1}^j(-1)^{k+1}\tbinom{2n+1}{j-k}\xi_k.
\end{equation}
 The Adams operations are determined by $\psi^t\xi_k=\xi_{kt}$.
\end{prop}

The next two theorems, \ref{algorithm} and \ref{comb}, determine
$v^{2m}(Spin(2n+1))$.
Two of the relations, $\car_2$ and $\car_3$, are defined in terms of the
 following algorithms. The algorithms will be justified in the proof of
 \ref{algorithm} and computed in \ref{comb}.
We use the relations $R_j$,
 $n+1\le j\le 2n$, and $S_j$, $2n<j\le 3n$, of
Corollary \ref{Spn}. Note that each equals $\pm\xi_j$ mod terms
with smaller subscripts of $\xi$.

For $\car_2$, begin with
\begin{equation}\label{b2}
\sum_{j=1}^n\sum_{k=1}^j(-1)^{k+1}\tbinom{2n+1}{j-k}\xi_{2k},\end{equation}
and subtract multiples of $R_{2n},\ldots,R_{n+1}$ to eliminate
$\xi_{2n},\ldots,\xi_{n+1}$. This results in a linear expression $E$ in
$\xi_1,\ldots,\xi_n$ which, we will show in (\ref{2pn}),
equals $(-1)^{n+1}2^n\xi_n$ plus lower
terms. Write (\ref{Ddef}) in the form
\begin{equation}\label{Ddef'}
2^{n+1}D+(-1)^n\xi_n+\sum_{k=1}^{n-1}(-1)^k\xi_k\sum_{i=0}^{n-k}\tbinom{2n+1}i.
\end{equation}
Add $2^n$ times (\ref{Ddef'}) to $E$ to eliminate $\xi_n$. The
resulting linear expression in $\xi_1,\ldots,\xi_{n-1}$, and $D$ is divisible
by $2^{n+1}$. Divide it by $2^{n+1}$, and then replace $\xi_{k}$ by $k^m$
(resp. 0) if $k$ is odd (resp. even) to get $\car_2$.

For $\car_3$, begin with
\begin{equation}\label{b3}
\sum_{j=1}^n\sum_{k=1}^j(-1)^{k+1}\tbinom{2n+1}{j-k}\xi_{3k},\end{equation}
and subtract multiples of $S_{3n},\ldots,S_{2n+1},R_{2n},\ldots,R_{n+1},$ and
of (\ref{Ddef'}) to
eliminate $\xi_{3n},\ldots,\xi_{n+1},\xi_n$. This results in
a linear expression in $\xi_1,\ldots,\xi_{n-1}$, and $D$ which is divisible by
$2^{n+1}$. Divide it by $2^{n+1}$ and replace $\xi_k$ by $k^m$ or 0 as
before to get an expression $E'$. The relation
$\car_3$ is $3^mD-E'$.

\begin{thm}\label{algorithm}
The abelian group $v^{2m}(Spin(2n+1))$ has generators $\xi_1$ and $D$ with
relations
\begin{equation} \label{Dreln}
2^{n+1}D=(\sum_{j=1}^n\sum_{\text{odd $k$}}\tbinom{2n+1}{j-k}k^m)\xi_1,
\end{equation}
$2^{eSp(m,n)}\xi_1$,
$\car_2$  and $\car_3$.

\end{thm}
\begin{pf} The second relation is the one obtained in Propositions \ref{vSpn}
and \ref{SUconj}.
This relation due to Adams operations on $\xi_i$-classes is still present in
$v^{2m}(Spin(2n+1))$. The first relation is the definition of $D$, together
with $\xi_k\sim\psi^k\xi_1\sim k^m\xi_1$ if $k$ is odd, and
$\xi_k\sim\psi^k\xi_1\sim0$ if $k$ is even.

The relations $\car_2$ and $\car_3$
are consequences of $\psi^2D\sim0$ and $\psi^3D\sim3^mD$, respectively. 
We evaluate
$\psi^2(2^{n+1}D)$ and $\psi^3(2^{n+1}D)$ using (\ref{Ddef}) and
$\psi^\ell\xi_k\sim\xi_{k\ell}$. The relations are reduced using $R_j$ and
$S_j$. Since $\psi^kD$ must be an integral combination of the generators
$\xi_1,\ldots,\xi_{n-1}$, and $D$, it follows that $\psi^2(2^{n+1}D)$ and
$\psi^3(2^{n+1}D)$, when expressed in terms of these classes, must be divisible
by $2^{n+1}$. Performing the division yields the desired relations for
$\psi^2(D)$ and $\psi^3(D)$.

By Corollary \ref{alls}, we need only to append the relation $\psi^{-1}(D)\sim
(-1)^mD$ in order to achieve all relations $\psi^k(D)\sim k^mD$ for all odd
$k$. However, by Proposition \ref{psi-1}, this relation adds no information if
$m$ is odd, while if $m$ is even, which will be treated at the end of this
section, it implies $2D\sim0$.\end{pf}

\begin{prop}\label{psi-1} $\psi^{-1}=-1$ in $PK^1(Spin(2n+1))$.\end{prop}
\begin{pf}
First note that $\psi^{-1}=(-1)^m$ in $PK^1(S^{2m+1})$. We use this to prove, by
 induction on
$n$, that $\psi^{-1}=-1$ in $PK^1(Sp(n))$.

Indeed, consider the short exact sequence
$$0\to PK^1(S^{4n-1})\mapright{p^*}PK^1(Sp(n))\mapright{i^*}PK^1(Sp(n-
1))\to0,$$
and assume $\psi^{-1}=-1$ in $PK^1(Sp(n-1))$. Let $y$ generate $\im(p^*)$.
Note $\psi^{-1}(y)=-y$. Let $x\in PK^1(Sp(n))$. By naturality of $\psi^{-1}$
and the induction hypothesis, we have $i^*(\psi^{-1}(x)+x)=0$, and hence
$\psi^{-1}(x)+x=\a y$ for some integer $\a$. Apply $\psi^{-1}$ and use
$\psi^{-1}\psi^{-1}=1$ to obtain $\a y=-\a y$, and hence $\a=0$.

Now we use the morphism $PK^1(Sp(n))\mapright{\theta} PK^1(Spin(2n+1))$
considered in Proposition \ref{Bmor} to deduce the
proposition. Since $2^{n+1}D\in\im(\theta)$, we obtain
$\psi^{-1}(2^{n+1}D)=-
2^{n+1}D$, and hence $\psi^{-1}(D)=-D$ since there is no torsion in
$PK^1(Spin(2n+1))$. All other elements are in the image from $PK^1(Sp(n)))$.
\end{pf}

The algorithm distills down
to the following closed form for the relations in Theorem \ref{algorithm},
which was observed from
extensive computer calculations by the second author and then proved (see
  Section \ref{combsec})
utilizing a number of elaborate combinatorial
arguments.
\begin{thm}\label{comb} $v^{2m}(Spin(2n+1))$ is the abelian group with
generators $\xi_1$ and
$D$ and four relations as below: $$2^{eSp(m,n)}\xi_1,$$
\begin{equation}\label{r1}\biggl(\sum_{\text{odd }k}k^m\sum_{i=0}^{n-
k}\tbinom{2n+1}i\biggr)\xi_1-2^{n+1}D,\end{equation}
\begin{equation}\label{r2}\biggl(\sum_{\text{odd }k}k^m\sum_{t\ge0}
\tbinom{2n+2}{n-1-k-4t}\biggr)\xi_1-2^nD,\end{equation}
and
\begin{equation}\label{r3}\biggl(2^n\sum_{\text{odd }k}k^m\sum_{t\ge0}
\tbinom{2n+1}{n-1-k-3t}\biggr)\xi_1-\tfrac13(2^{2n+1}+1-3^{m+1})D.
\end{equation}
If $m$ is odd and $n\ne4$, the coefficients of $\xi_1$ in $(\ref{r1})$ and
$(\ref{r2})$ are divisible by $2^n$.
\end{thm}

From this, we easily deduce Theorem \ref{main}. Subtract 2 times (\ref{r2})
from (\ref{r1}) to obtain the relation $R_1(m,n)$ of \ref{main} on the
generator $\xi_1$. Note that the exponent of 2 in the coefficient of $D$ in
(\ref{r3}) is $\min(2n+1,\nu(m+1)+2)$ if $m$ is odd. If $n\le\nu(m+1)+2$, 
divide (\ref{r2}) by $2^n$ 
to split off $\bz/2^n$, the second summand in \ref{main}, and then
add an appropriate multiple of (\ref{r2}) to (\ref{r3}) to eliminate $D$,
yielding the relation $R_2(m,n)$ of \ref{main} on $\xi_1$. The last part of
Theorem \ref{comb} is necessary in order to know that we can split off
$\bz/2^n$.

Similarly, if $\nu(m+1)+2<n$, use (\ref{r3}) to split off $\bz/2^{\nu(m+1)+2}$,
and then subtract an appropriate multiple of (\ref{r3}) from (\ref{r2}) to
obtain $R_2(m,n)$. This completes the proof of Theorem \ref{main} when $m$ is
odd.

To prove that $v^{2m}(Spin(2n+1))\approx\bz/2\oplus\bz/2$ when $m$ is even, we
begin by showing that $eSp'(m,n)=1$ in Corollary \ref{vSpn} when $m$ is even.
It suffices to show that $R_{m,2n}\equiv2$ mod 4 and that $R_{m,j}$ and
$S'_{m,j}$ are even for all relevant $j$. Since $k^m\equiv1$ mod 4 if $k$ is odd
and $m$ even, we have, mod 4,
$$S'_{m,j}\equiv\sum_{k\text{ odd}}\tbinom jk=2^{j-1},$$
$$R_{m,2n}\equiv\sum_{k\text{ odd}}\tbinom{2n+1}k-2\tbinom{2n+1}{2n+1}=2^{2n}-
2,$$
and other $R_{m,j}$ are congruent to $2^{2n}$ minus 2 times a sum of binomial
coefficients.

As observed in the proof of \ref{algorithm},
$\psi^{-1}$ implies $2D\sim0$ when $m$ is even. Thus it
remains to show that the coefficients of relations in Theorem \ref{comb} are
even when $m$ is even, which reduces to showing
\begin{equation}\label{e1}\sum_{k\text{ odd}}\sum_{i=0}^{n-
k}\tbinom{2n+1}i\equiv0\mod 2\end{equation}
and
\begin{equation}\label{e2}\sum_{k\text{ odd}}\sum_{t\ge0}\tbinom{2n+2}{n-
k-1-4t}\equiv0\mod 2.\end{equation}
The sum in (\ref{e1})
is $\sum_{i=0}^{n-1}\tbinom{2n+1}i\bigl[\frac{n-i+1}2\bigr]$. If
$n$ is even, terms $2j$ and $2j+1$ in this sum will have the same parity.
If $n\equiv3$ mod 4, terms $2j+1$ and $2j+2$ will have the same parity.
If $n\equiv1$ mod 4, the only odd terms will occur in pairs $i=8j$ and
$i=8j+3$.

The sum in (\ref{e2}) is $\sum\tbinom{2n+2}{n-2j}\bigl[\frac{j+1}2\bigr]$.
If $n$ is odd, all coefficients $\binom{2n+2}{n-2j}$ are even, while if $n=2N$,
this becomes
$\sum\tbinom{2N+1}{N-j}\bigl[\frac{j+1}2\bigr]$ mod 2,
which equals the first sum,
and hence is even. This completes the proof of Theorem \ref{main}.

\section{Proof of Theorem \ref{comb}}\label{combsec}
In this section, we perform the algorithm in Theorem \ref{algorithm} and obtain
Theorem \ref{comb}. First note that (\ref{r1}) is an elementary manipulation of
(\ref{Dreln}).

We begin by deriving (\ref{r2}). Let $M$ be a matrix whose rows present the
relations $R_{2n},\ldots,R_{n+1}$ of \ref{Spn} expressed in terms of
$\xi_{2n},\ldots,\xi_1$. Then
$$M_{i,j}=(-1)^{j-i}\tbinom{2n+1}{j-i}+(-1)^{i+j}\tbinom{2n+1}{j+i-2n-1}.$$
Let $M_L$ (resp. $M_R$) denote the left (resp. right) half of $M$.

Let $\sigma_j:=\sum_{k=0}^j\binom{2n+1}k$. Let
$P$ denote a row vector of length $2n$ whose $(2j+1)$st entry is
$(-1)^{n+1+j}\sigma_j$ for $0\le j<n$, with other entries 0. This represents
(\ref{b2}). Write $P=(P_L|P_R)$, with $P_L$ and $P_R$ of length $n$.
We wish to add multiples of the rows of $M$ to $P$ to annihilate $P_L$.
This is facilitated by first replacing $M$ by
the equivalent set of relations $M_L^{-1}M=[I|M_L^{-1}M_R]$.
The linear expression $E$ in the algorithm for $\car_2$ described
before Theorem \ref{algorithm}
corresponds to $P_R-P_LM_L^{-1}M_R$. We will show that this equals
\begin{equation}\label{vec}(-1)^{n+1}2^n\sum_{j=0}^{n-1}(-
1)^j\tbinom{2n+1}jV_j,\end{equation}
where $V_j$ is a vector of length $n$ which begins with $j$ 0's followed by
$$1,1,-1,-1,1,1,-1,\ldots.$$

The expression corresponding to (\ref{vec}) is
\begin{equation}\label{2pn}
2^n\bigl(\sum_{k=1}^n(-1)^{k+1}\xi_k\bigl(\sum_{j=0}^{n-k}\tbinom{2n+1}j
-2\sum_{t\ge0}\tbinom{2n+2}{n-1-k-4t}\bigr)\bigr).\end{equation}
We have used
\begin{equation}\label{Pas}\tbinom{2n+2}a=\tbinom{2n+1}a+\tbinom{2n+1}{a-1}
\end{equation}
to show that this corresponds to (\ref{vec}).
Adding $2^n$ times (\ref{Ddef'}) to (\ref{2pn}) and replacing $\xi_k$ by $k^m$
or 0 yields $-2^{n+1}$ times (\ref{r2}), as
desired.

To complete the proof of (\ref{r2}), it remains to prove (\ref{vec}).
Let $c_i=\binom{2n+1}i$ and $d_i=\binom{2n+i}i$. Consideration of generating
functions implies that for any positive integer $m$
\begin{equation}\label{zero}\sum_{j=0}^m(-1)^jd_jc_{m-j}=0\end{equation}
and that the $(i,j)$th entry of $M_L^{-1}$ is $d_{j-i}$. Using these facts, one
can easily prove the following lemma.
\begin{lem}\label{DCsum}$M_L^{-1}M_R=\sum_{j=1}^nD_jC_j$, where $D_j$ is a
column vector of length $n$ with $t^{\text{th}}$
entry $d_{n+1-t-j}-d_{n+j-t}$, and $C_j$
is a row vector of length $n$ with $t^{\text{th}}$ entry $(-1)^{t-j}c_{t-j}$.
\end{lem}
To simplify notation, we shall complete the proof of (\ref{vec})
when $n$ is odd. The proof
when $n$ is even is virtually the same. Let $n=2e+1$. We shall prove the
following lemma.
\begin{lem}\label{PD} For $1\le j\le n$,
$$-P_LD_j=(-1)^{[(j-1)/2]}2^n+\sum_{u=1}^{[j/2]}(-1)^{e+u+1}\sigma_{e+u}d_{j-
2u}.$$
\end{lem}
Then the entry in the $t^{\text{th}}$ position of $-P_LM_L^{-1}M_R$ is
$$(-P_L\sum_jD_jC_j)_t=A_{1,t}+A_{2,t},$$
where
$$A_{1,t}=2^n\sum_{j=1}^t(-1)^{[(j-1)/2]}(-1)^{t-j}c_{t-j},$$
which equals the $t^{\text th}$ entry of (\ref{vec}), and
$$A_{2,t}=\sum_{j=1}^t(-1)^{t-j}c_{t-j}\sum_{u=1}^{[j/2]}(-
1)^{e+u+1}\sigma_{e+u}d_{j-2u}.$$
Using (\ref{zero}), one can check that
$$A_{2,t}=\begin{cases}0&t\text{ odd}\\(-1)^{e+T+1}\sigma_{e+T}&t=2T.
\end{cases}$$
Thus $P_R+A_2=0$, and hence $P_R-P_LM_L^{-1}M_R$ equals (\ref{vec}), as
claimed. Once we have given the proof of Lemma \ref{PD}, we will be done with
the proof of the $\car_2$-part of Theorem \ref{comb}.

\begin{pf*}{Proof of Lemma \ref{PD}} We use generating functions.
The coefficients of $x^0,\ldots,x^{n-1}$ in $p(x):=(1-x^2)^{2n+1}/(1+x^2)$ give
the
entries of $P_L$, while $q(x):=(1-x)^{-(2n+1)}=\sum d_ix^i$.
Note that $p(x)q(x)=(1+x)^{2n+1}/(1+x^2)$. We obtain
\begin{eqnarray}-P_LD_j&=&\sum_{t=1}^n\coef(x^{t-1},
p(x))\cdot\bigl(\coef(x^{n+j-t},q(x))-\coef(x^{n-j+1-t},q(x)
\bigr)\nonumber\\
&=&\coef(x^{n+j-1},p(x)q(x))-\coef(x^{n-j},p(x)q(x))
\label{coefs}\\
&&-\sum_{t>n}\coef(x^{t-1},p(x))\cdot d_{n+j-t}.\nonumber
\end{eqnarray}
The sum subtracted at the end is due to the truncation of $D_j$; one easily
checks that this equals the $u$-sum in Lemma \ref{PD}. Finally we observe that,
if
we employ the symmetry property of binomial coefficients, then the expression
in line (\ref{coefs}) equals the following expression if $j\equiv 1$ or 2 mod 4, and
equals the negative of this expression otherwise.
$$\tbinom{2n+1}n-\tbinom{2n+1}{n-1}-\tbinom{2n+1}{n-2}+\tbinom{2n+1}{n-
3}+\tbinom{2n+1}{n-4}-\tbinom{2n+1}{n-5}-\cdots$$
This expression equals $2^n$, as can be seen by expanding $(1-i)^{2n+1}$ or by
induction on $n$.\end{pf*}

The derivation of (\ref{r3}) is similar to that of (\ref{r2}) just performed,
although somewhat more complicated. Replacing the matrix $M$ used in (\ref{r2})
is a
matrix $N$ whose rows present the relations
$S_{3n},\ldots,S_{2n+1},R_{2n},\ldots,R_{n+1}$ of \ref{Spn} in terms of
$\xi_{3n},\ldots,\xi_1$. Then
$$N_{i,j}=\begin{cases}(-1)^{j-i}\tbinom{3n+1-i}{j-i}&1\le i\le n\\
(-1)^{j-i}\tbinom{2n+1}{j-i}+(-1)^{i+j}\tbinom{2n+1}{j+i-4n-1}&n<i\le 2n.
\end{cases}$$
Divide $N$ into $n\times n$ submatrices, which we name as follows:
$$N=\begin{pmatrix} T_1&U_1&U_2\\0&T_2&U_3\end{pmatrix}.$$
We perform row operations on $N$ to reduce it to
$$N'=\begin{pmatrix}I&0&B_2-B_1W\\
0&I&W\end{pmatrix},$$
where $B_i=T_1^{-1}U_i$, and $W=
T_2^{-1}U_3=\sum_{j=1}^n D_jC_j$, exactly as in Lemma \ref{DCsum}.

Let $\sigma_j,c_i,d_i,C_j$, and $D_j$ be as in the preceding proof,
with $\si_j=0$ if $j$ is not an integer.
Let $Q$ denote a row vector of length $3n$ whose $(3j+1)$st entry is
$(-1)^{n+1+j}\sigma_j$ for $0\le j<n$, with other entries 0. This represents
(\ref{b3}). Write $Q=(Q_0|Q_1|Q_2)$ with each $Q_i$ of length $n$.
We use $N'$ to reduce $Q$ to $(0|0|Q'_2)$, where
\begin{equation}\label{Q2'}Q_2'=Q_2-Q_0(B_2-B_1W)-Q_1W.\end{equation}
We will show this equals
\begin{multline}\label{Q2"}
(-1)^n\biggl({\tfrac13}(2^{2n+1}+1)(-\si_0,\si_1,-\si_2,\si_3,\ldots,\si_{n-
1})\\
+2^{2n+1}
(0,-c_0,c_1,-c_2,c_3+c_0,-(c_4+c_1),\ldots)\biggr).\end{multline}

Thus the vector (\ref{Q2"}) expresses $\psi^3(2^{n+1}D)$ in terms of
$\xi_n,\ldots,\xi_1$ as
\begin{equation}\label{two}{\tfrac13}(2^{2n+1}+1)\sum_{k=1}^n(-
1)^{k+1}\xi_k\sum_{i=0}^{n-k}\tbinom{2n+1}i+2^{2n+1}\sum_{k=1}^{n-1}
(-1)^k\xi_k\sum_{t\ge0}\tbinom{2n+1}{n-1-k-3t}.\end{equation}
We add $\frac13(2^{2n+1}+1)$ times (\ref{Ddef'}) to eliminate $\xi_n$ and get
$${\tfrac13}(2^{2n+1}+1)2^{n+1}D+2^{2n+1}\sum_{k=1}^{n-1}
(-1)^k\xi_k\sum_{t\ge0}\tbinom{2n+1}{n-1-k-3t}.$$
Divide by $2^{n+1}$, replace $\xi_k$ by $k^m$ or 0, and subtract from $3^mD$;
this yields exactly (\ref{r3}).

It remains to show that $Q_2'$ equals (\ref{Q2"}). We show this when $n\equiv0$
mod 6; other congruences can be handled similarly. Let $r=n/3$.
It is not difficult to prove that
$$(B_1)_{i,j}=(-1)^j\tbinom{n+j-i-1}{n-i}\tbinom{3n+1-i}{n+j-i}\text{ and }
(B_2)_{i,j}=(-1)^j\tbinom{2n+j-i-1}{n-i}\tbinom{3n+1-i}{2n+j-i}.$$
We obtain that, for $0\le t<n$, the $(t+1)$st entry of $Q_2'$ is
\begin{eqnarray}\label{step0}&&(-1)^{t+1}\si_{2r+t/3}-\sum_{u=0}^{r-1}
(-1)^{u+t}\si_u\tbinom{2n+t-3u-1}{n+t}\tbinom{3n-3u}{n-t}\\
\label{triple} &+&\sum_{s=0}^{n-1}\sum_{u=0}^{r-1}(-1)^{u+s}
\sigma_u\tbinom{n+s-3u-
1}{s}\tbinom{3n-3u}{2n-s}\\
\nonumber&&\qquad\qquad\qquad\qquad
\cdot\sum_{j=0}^{t}(d_{n-s+j-t-1}-d_{n+t-j-s})(-1)^{j}c_{j}\\
\label{double} &+&\sum_{u=0}^{r-1}(-1)^u\si_{r+u}\sum_{j=0}^t
(d_{n-3u-1+j-t}-d_{n+t-j-3u})(-1)^jc_j.
\end{eqnarray}
It follows from (\ref{zero}) that the $j$-sum in (\ref{triple}) is 0  for
$s\ge n$, except that it is $-1$ when $s=n+t$. Hence (\ref{triple}) can be made
into a sum over all $s\ge0$ with the sum-part
of (\ref{step0}) incorporated into it.

Next we note that if $u=r+\delta$ with $\delta\ge0$ in (\ref{double}), then,
again using (\ref{zero}), the $j$-sum in (\ref{double}) is nonzero only when
$\delta=t/3$, in which case it is $-1$. Thus
(\ref{double}) can be extended to
$u\ge0$ if the initial term of (\ref{step0}) is included.
Similarly we note that if $u=r+\Delta$ with $\Delta\ge0$ in (\ref{triple}), then
$\tbinom{n+s-3u-1}{s}\tbinom{3n-3u}{2n-s}$ is nonzero only when $s=3\Delta$, in
which case it is $(-1)^{\Delta}$. Thus the modified
(\ref{double}) can be incorporated into
(\ref{triple}) by extending the $u$-sum
of (\ref{triple}) to $u\ge0$. Now the $(t+1)$st entry of
$Q_2'$ has been simplified to
\begin{equation}\label{trip}
\sum_{s\ge0}\sum_{u\ge0}(-1)^{u+s}\si_u\tbinom{n+s-3u-
1}{s}\tbinom{3n-3u}{2n-s}\sum_{j=0}^{t}(d_{n-s+j-t-1}-d_{n+t-j-s})(-1)^{j}c_{j}.
\end{equation}

We will prove the following lemma at the end of this section.
\begin{lem} \label{binomlem} If $n\ge0$, and $v$ and $w$ are arbitrary
integers, then
$$\sum_{s\ge0}(-1)^s\tbinom{n+s-1-w}s\tbinom{3n-w}{2n-s}\tbinom{3n-s+v}{2n}=
\tbinom{4n-w+v}{2n}.$$
\end{lem}
Using this with $w=3u$, (\ref{trip}) reduces to
\begin{equation}\label{trip2} \sum_{j=0}^t (-1)^jc_j\sum_{u\ge0}(-1)^u\si_u
(\tbinom{4n-3u+j-t-1}{2n}-\tbinom{4n-3u+t-j}{2n}).\end{equation}
Viewing $\binom{4n-3u+\eps}{2n}$ as $(-1)^{u+\eps}\coef(x^{2n-3u+\eps},(1+x)^{-
(2n+1)})$, and $\si_u$ as $$\coef(x^{3u},(1+x^3)^{2n+1}/(1-x^3)),$$ we obtain a
simplification of (\ref{trip2}) to
\begin{equation}\label{trip3}
(-1)^{t+1}\sum_{j=0}^t c_j\bigl(\coef(x^{2n+j-t-1},\tfrac{(1-x+x^2)^{2n+1}}
{1-x^3})+\coef(x^{2n+t-j},\tfrac{(1-x+x^2)^{2n+1}}{1-x^3})\bigr).\end{equation}

We will prove the following lemma at the end of this section.
\begin{lem}\label{sym} Let $f_m(x)=(1-x+x^2)^{m+1}/(1-x^3)$. Then for
$m\ge0$ and $j\ge0$,
$$\coef(x^{m-j-1},f_m(x))+\coef(x^{m+j},f_m(x))=\begin{cases} (1+2(-2)^m)/3
&j\equiv0,2\mod 3\\ (1+2(-2)^{m+1})/3&j\equiv1\mod 3.\end{cases}$$
\end{lem}
We obtain that (\ref{trip3}) equals
$$(-1)^{t+1}\biggl({\tfrac13}(2^{2n+1}+1)\sum_{j=0}^t c_j-2^{2n+1}\sum
\begin{Sb}0\le j\le t\\t-j\equiv1\ (3)\end{Sb} c_j\biggr),$$
which equals the component of (\ref{Q2"}) in position $t+1$, as desired,
completing the proof of (\ref{r3}), as described surrounding (\ref{Q2"}).

Next we prove the portion of Theorem \ref{comb} which states that
\begin{equation}\label{summ}
\sum_{\text{odd }k}k^m\sum_{t\ge0}\tbinom{2n+2}{n-1-k-
4t}\end{equation}
is divisible by $2^n$ if $m$ is odd and $n\ne4$. This occurs in (\ref{r2}).
We accomplish this by showing 
later in this section that, if $m=2a+1$, then (\ref{summ}) equals
\begin{equation}\label{bigsum}
\tbinom{n-1}12^{2n-4}+\sum_{j=2}^{n/2}\tbinom{n-j}j2^{2n-4j}\sum_{i\ge j-1}
8^i\tbinom ai\sum_{t=0}^{j-2}(-1)^t\tbinom{2j-1}t(2j-2t-1){\tbinom{j-t}2}^i.
\end{equation}
Then we note that the first term of (\ref{bigsum}) is divisible by $2^n$
provided $n\ge3$. For each value of $j\ge2$, the $j$-summand is nonzero for
$n\ge 2j$ and is clearly divisible by $2^{2n-4j+3(j-1)}$, which is $\ge2^n$
provided $n\ge j+3$. The case $n=4$ and $j=2$ is the only time that this gives
a nonzero term which might not be divisible by $2^n$; it is what causes the
restriction $n\ne4$ in the last part of Theorem \ref{comb}.

The expression (\ref{bigsum}), while not as elegant as (\ref{summ}), is often
more useful in performing specific calculations of $\nu(-)$, such as in
Proposition \ref{5and6}. Before proving (\ref{bigsum}), we handle the case
$n=4$, which was excluded from Theorem \ref{main}.

\begin{prop}\label{n=4}
$$v^{2m}(Spin(9))\approx\begin{cases}\bz/2\oplus\bz/2&m\text{ even}\\
\bz/8\oplus\bz/2^{e(m)}&m\text{ odd},\end{cases}$$
where
$$e(m)=\begin{cases}\min(\nu(m-5)+2,6)&\text{if $m\equiv1$ mod $4$}\\
\min(\nu(m-7)+2,8)&\text{if $m\equiv3$ mod $4$}.\end{cases}$$
\end{prop}
\begin{pf} The numbers $eSp(m,n)$ are as in \cite[5.1]{BDM} (suitably
interpreted), and large enough as to not matter here. Let $P=3^m-3$.
The other three relations of \ref{comb} are $(160+10P)\xi_1-32D$, $(48+P)\xi_1-
16D$, and $16(39+P)\xi_1+(P-168)D$. If $m\equiv3$ mod 4, then $P\equiv8$
mod 16, so we get a $\bz/8$ due to $(48+P)\xi_1$. (This is the anomaly.) If
$m\equiv1$ mod 4, $P\equiv0$ mod 16, and we get a $\bz/8$ due to $(P-168)D$.

If $m\equiv 1$ mod 4, let $Q=3(3^{m-5}-1)$, so that $\nu(Q)=\nu(m-5)+2$.
Replace $P$ by $3^4Q+3^5-3$.
Subtract multiples of the third relation from the other two so as to eliminate
the $D$-term from each. The relations obtained are, up to odd multiples,
$2^6+2Q+Q^2/4$ and $2^8+Q+Q^2/8$. For $\nu(Q)\ge4$, the minimum of the exponent
of 2 of these two expressions is $\min(6,\nu(Q))$, as claimed. The case
$m\equiv3$ mod 4 is handled similarly.\end{pf}

Now we establish the equivalence of (\ref{summ}) and (\ref{bigsum}). A key
ingredient is the following lemma.
\begin{lem} \label{Afor} If $n\ge0$ and $A$ is an integer, then
$$\sum_{j=0}^{n/2}(-1)^j\tbinom{n-j}j4^{n-2j}\tbinom{2j-1}{j-A}=(-
1)^A\sum_{t\ge0}\tbinom{2n+2}{n-2A-4t}.$$
\end{lem}
\begin{pf} Thinking of either expression as a function $f(n,A)$, we establish
that each satisfies
$$f(0,A)=\begin{cases}0&A>0\\ 1&A\le0,\text{ even}\\
-1&A<0,\text{ odd,}\end{cases}$$ $f(-1,A)=0$,
and the recursive formula
$$f(n,A)=4f(n-1,A)-f(n-2,A-1)-2f(n-2,A)-f(n-2,A+1).$$
Thus they are equal. The recursive formula for the LHS of \ref{Afor}
is proved using several
applications of (\ref{Pas}), while for the RHS one proves the closely related
formula
$$\tbinom{2n+2}m=4\tbinom{2n}{m-1}+\tbinom{2n-2}m-2\tbinom{2n-2}{m-
2}+\tbinom{2n-2}{m-4}$$
by multiplying both sides by $x^m$ and summing over all values of $m$.
\end{pf}

Replace $k$ in (\ref{summ}) by $2s+1$, and $m$ by $2a+1$.
Using \ref{Afor} with $A=s+1$,
(\ref{summ}) becomes
\begin{eqnarray*}&&\sum_{s\ge0}(2s+1)^{2a+1}(-1)^{s+1}\sum_{j=0}^{n/2}
(-1)^j\tbinom{n-j}j
4^{n-2j}\tbinom{2j-1}{j-s-1}\\&=&\sum_{j=0}^{n/2}(-1)^j\tbinom{n-j}j4^{n-
2j}\sum_{s\ge0}(-1)^{s+1}(2s+1)\tbinom{2j-1}{j-s-1}\sum_{i\ge0}\tbinom
ai(4s^2+4s)^i.
\end{eqnarray*}
Replacing $s$ by $j-t-1$ yields a formula which can easily be manipulated to
(\ref{bigsum}), except that the $i$-sum is unrestricted. The restriction to
$i\ge j-1$ follows from the following lemma.
\begin{lem} If $1\le i\le j-2$, then
$$\sum_{t=0}^{j-2}(-1)^t\tbinom{2j-1}t(2j-2t-1){\tbinom{j-t}2}^i=0.$$
\end{lem}
\begin{pf}We solve the system of equations
\begin{eqnarray*}a_1\tbinom32\ +\ a_2\tbinom 42\quad+\cdots+\ a_{j-2}\tbinom j2
\quad&=&1\\
a_1{\tbinom32}^2\ +\ a_2{\tbinom 42}^2\
+\cdots\ +\ a_{j-2}{\tbinom j2}^2\ &=&1\\
&\vdots&\\
a_1{\tbinom 32}^{j-2}+a_2{\tbinom42}^{j-2}+\cdots+a_{j-2}{\tbinom j2}^{j-2}&=&1
\end{eqnarray*}
by Cramer's rule and the Vandermonde evaluation of determinants, obtaining
\begin{eqnarray*}a_{s-2}&=&\frac1{\tbinom s2}\frac{\prod_{t>s}\bigl(\tbinom t2-
1\bigr)\prod_{t<s}\bigl(1-\tbinom t2\bigr)}{\prod_{t>s}\bigl(\tbinom t2-
\tbinom s2\bigr)\prod_{t<s}\bigl(\tbinom s2-\tbinom t2)\bigr)}\\
&=&\frac2{s(s-1)}\frac{\prod_{t\ne s}(t(t-1)-2)}{\prod_{t\ne s}(t(t-1)-s(s-
1))}\\
&=&\frac2{s(s-1)}\frac{\prod_{t\ne s}(t-2)(t+1)}{\prod_{t\ne
s}(t-s)(t+s-1)}\\
&=&(-1)^{s+1}\frac{2s-1}3\frac{\tbinom{2j-1}{j-s}}{\tbinom{2j-1}{j-2}}.
\end{eqnarray*}
The values of $t$ in the products always run from 3 to $j$. Substituting the
solutions into the equations yields equations which reduce to those that were
to be proved.\end{pf}

This completes the proof that the coefficient of $\xi_1$ in (\ref{r2}) is
divisible by $2^n$ if $n\ne 4$. The proof that the coefficient of $\xi_1$ in
(\ref{r1}) is also divisible by $2^n$
is extremely similar. Using virtually identical methods, we show it
equals
$$\tbinom{n+1}12^{2n-3}+\sum_{j\ge2}2^{2n+1-4j}\biggl(\tbinom{n+2-j}j-\tbinom
{n-j}{j-2}\biggr)\sum_{i\ge j-1}8^i\tbinom ai\sum_{t=0}^{j-2}(-1)^t\tbinom{2j-
1}t(2j-2t-1){\tbinom{j-t}2}^i.$$
Each term of this is divisible by $2^{n+1}$ except when $n=4$ and $j=2$, in
which case it is divisible by $2^n$.

We complete the paper by giving deferred proofs of two lemmas.
\begin{pf*} {Proof of Lemma \ref{binomlem}} This lemma can probably be proved
by some sort of generating function argument. However, we present a proof along
the lines by which it was discovered, utilizing {\tt Maple} and the method of
\cite{WZ}.

As described in \cite{WZ}, one can often prove $f(n):=\sum_s F(n,s)$ is constant
(independent of $n$) by finding a function $G$ such that
\begin{equation}\label{recur}F(n+1,s)-F(n,s)=G(n,s+1)-G(n,s).\end{equation}
If such a function $G$ can be found, then summing this equation over all values
of $s$ yields the conclusion $f(n+1)-f(n)=0$, and hence $f$ is constant.

In our application, we run the software associated with the book \cite{WZ}
on the function
\begin{equation}F(n,s):=(-1)^s\tbinom{n+s-1-w}s\tbinom{3n-w}{2n-s}\tbinom{3n-
s+v}{2n}\biggm/\tbinom{4n-w+v}{2n}.\label{Fns}\end{equation}
Here $n$, $w$, and $v$ are all free parameters, but we tell {\tt Maple} that
$n$ is the induction parameter by running
{\tt zeil(F(n,s),s,n,N)}, where $F$ is as in (\ref{Fns}). The second argument
of {\tt zeil} is the summation variable, the third is the parameter to be used
for the recursion, and the fourth is the symbol to use for the operator which
increases the value of the third parameter by 1. The output in this case
is {\tt N-1,} followed by $G(n,s)$, an explicit polynomial, involving $v$ and
$w$ in addition to $n$ and $s$, for which (\ref{recur}) holds.  In this case,
the polynomial $G$ involves approximately 500 monomials. The significance of
the {\tt N-1} is that this operator acting on $F$ equals $G(n,s+1)-G(n,s)$.
This says exactly (\ref{recur}). So we know that the ratio of the two
expressions involved in the lemma is constant. To know that this ratio is 1,
we note that when $n=0$ each equals 1.\end{pf*}

\begin{pf*} {Proof of Lemma \ref{sym}} The proof is by induction on $m$.
The lemma is readily verified if $m=0$. Let $s_{m,j}$ denote the sum of the two
coefficients with which the lemma deals. Since $f_{m+1}=(1-x+x^2)f_m$, we
obtain $$s_{m+1,j}=\begin{cases}s_{m,1}&\text{if $j=0$}\\
s_{m,j-1}-s_{m,j}+s_{m,j+1}&\text{if $j>0$.}\end{cases}$$
The induction step follows easily.\end{pf*}

\def\line{\rule{.6in}{.6pt}}

\end{document}